\documentclass{article}

\usepackage{amsmath}
\usepackage[francais]{babel}
\usepackage{amssymb}
\usepackage[applemac]{inputenc}

\def\n{\noindent}
\def\Q{\hbox{\bf Q}}

\def\F{\hbox{\bf F}}

\def\Z{\hbox{\bf Z}}
\def\C{\hbox{\bf C}}
\def\G{\Gamma}

\def\sup{\mathop{\rm sup}}
\def\lim{\mathop{\rm lim}}

\def\<{\mathop{<\!\!<}}

\def\GL{\mathop{\bf GL}}

\def\Gal{\mathop{\rm Gal}}

\def\Aut{\mathop{\rm Aut}}

\def\ind{\rm{ind}}
\def\l{\ell}
\def\r{\rho}

\def\to{\rightarrow}

\def\g{\gamma}
\def\n{\noindent}

\def\kbar{\overline{k}}
\def\vbar{\overline{v}}
\def\Xbar{\overline{X}}

\def\se{\mbox{\bf--}}

\begin{document}

\title { Un critère d'indépendance pour une famille de représentations $\l$-adiques}

\author {Jean-Pierre Serre}

\maketitle

   \bigskip
   
   \date
   
      {\bf Introduction}
      
      \bigskip

      Soit  $k$  un corps de nombres, de clôture algébrique $\kbar$, et soit $A$ une variété abélienne sur  
     $k$, de dimension $d$. Comme on sait, de telles données définissent, pour tout nombre premier $\l$, une 
     représentation $\l$-adique
     
     \smallskip
     
     $ \r_\l : \G_k \ \to \ \Aut(T_\l(A)) \ \cong \ \GL _{2d}(\Z_\l),$
     
     \smallskip
     
     \n où  $\G_k = \Gal(\kbar/k)$, et $T_\l(A)$ est le $\l$-ième module de Tate de $A$ sur $\kbar$. La famille des   $ \r_\l $ s'identifie à un homomorphisme continu
     
      \smallskip
     $ \r  : \  \G_k \ \to \ \prod_\l \Aut(T_\l(A)) \ \cong \ \ \prod_\l \GL _{2d}  (\Z_\l)$.
     
      \smallskip
      
    Lorsqu'on s'intéresse au sous-groupe $\r(\G_k)$ de $\prod_\l \GL _{2d}(\Z_\l)$, il est commode de savoir que  $\r(\G_k)$ est le produit direct des
     $\r_\l(\G_k)$, autrement dit, que les $\r_\l$ sont ``indépendants". Bien entendu ce n'est pas toujours vrai, mais on peut démontrer (cf. [Se 86]) que cela le devient après une extension finie convenable de $k$; autrement dit,  les $\r_\l$ sont ``presque indépendants".
     
     Je me propose de reprendre cette question en mettant en évidence les propriétés des $\r_\l$ qui entraînent la presque indépendance. Comme on le verra au \S2,  ce sont des propriétés de ramification, analogues à ce que l'on appelle la ``semi-stabilité''; curieusement, les éléments de Frobenius, si utiles en d'autres circonstances, ne jouent ici aucun rôle.

     L'intérêt de cette axiomatisation est qu'on peut l'appliquer à des situations plus générales que celle des variétés abéliennes, par exemple à la cohomologie $\l$-adique des variétés algébriques sur un corps de nombres, cf. \S 3.2.  Un résultat très voisin avait d'ailleurs été obtenu il y a une quinzaine d'années par M.J. Larsen et R. Pink dans des lettres (datées du 23/5/95 et 26/5/95) dont le contenu n'a malheureusement pas été publié jusqu'à présent.
     
     \bigskip
        La démonstration du théorème principal (théorème 1 du \S2) est donnée 
     au \S8. Elle repose sur diverses propriétés des corps de nombres et des groupes linéaires (corps de classes, théorème de Hermite-Minkowski,
     théorèmes de Jordan et de Nori); ces propriétés font l'objet des \S\S4-7.
     
     \smallskip
     
 \n    {\it Remerciement.} Cet article doit beaucoup à L. Illusie: il m'a encouragé à l'écrire, il m'a fourni de nombreuses références et il m'a communiqué ([Il 10]) une
démonstration d'un résultat auxiliaire essentiel, qui avait été démontré auparavant, sous une forme un peu différente, par 
N. Katz et G. Laumon. Je lui en suis très reconnaissant.

   \bigskip
  \n {\bf \S1. La notion d'indépendance.}
   
   \smallskip
   
 % [ On travaille dans la catégorie des groupes profinis; les homomorphismes
 %  sont supposés continus et les sous-groupes sont supposés fermés.]
   
     Soit  $\G$ un groupe, et soit  $\r_i : \G \to G_i$ une famille d'homomorphismes de $\G$ dans des groupes $G_i$ indexés par un ensemble $I$. Cela revient à se donner un homomorphisme 
     
 \quad \quad     $\r = (\r_i) : \G  \ \to  \ {\small{\prod}}_{i \in I} G_i .$ 
 
 \smallskip
      
   \n   On dit que les $\r_i$ sont {\it indépendants} si la propriété suivante est
satisfaite:

\smallskip

     (R) \quad $\r(\G) = \prod \r_i(\G).$
 \smallskip
 
\n Autrement dit, si  $\g_i$  est une famille quelconque d'éléments de $\G$, il existe  $\g \in \G$  tel que  $\r_i(\g)= \r_i(\g_i)$  pour tout  $i$.
 
 \smallskip
 
 Il y a une propriété plus faible que l'on peut considérer:

\smallskip

  (RO)  \quad $\r(\G)$  est un sous-groupe d'indice fini de $ \prod \r_i(\G)$.
  
  \smallskip
  
   A partir de maintenant,  on suppose que $\G$ est un groupe profini, que les  $G_i$ sont localement
   compacts, et que les  $\r_i$  sont continus (de sorte que les  $\r_i(\G)$
   sont des groupes profinis). On s'intéresse à la propriété:
  
  \smallskip
  
  (PR) Il existe un sous-groupe ouvert  $\G'$ de  $\G$  tel que les restrictions
  des  $\r_i$  à  $\G'$  vérifient (R). [Noter que  $\G'$ est d'indice fini dans  $\G$, puisque $\G$ est compact.] 
  
  On dit alors que les $\r_i$ sont {\it presque indépendants}.
 
 \smallskip
 
   On a  (R) $\Rightarrow$ (RO) $\Rightarrow$ (PR): c'est clair pour  (R) $\Rightarrow$ (RO), et
   ce n'est pas difficile pour  (RO) $\Rightarrow$ (PR).
   
   \smallskip
   
   {\it Remarque.} On peut aussi exprimer  (R)  comme une propriété des noyaux  $N_i$  des $\r_i$. Si l'on pose  $N'_i = \bigcap_{j\neq i}N_j$, la
   condition  (R)  est équivalente à chacune des deux conditions
   suivantes:
   
    (R1)  $\G = N_i .N'_i$  pour tout  $i$.
    
    (R2)  $\G$ est engendré (topologiquement) par les $N'_i$.
    
    \n[Lorsque  $I$  est fini, cela se démontre par récurrence sur le nombre 
    d'éléments de $I$; le cas général s'en déduit par passage à la
    limite, en utilisant la compacité de  $\G$.]
    
      On peut préciser  (R2): si l'on note  $\G'$ le plus petit sous-groupe
      fermé de  $\G$  contenant les $N'_i$, alors  $\G'$  est le plus grand
      sous-groupe fermé de $\G$ sur lequel les  $\r_i$  sont indépendants.

   \bigskip
   
  \n {\bf \S2. Enoncé du théorème.}
   
   \smallskip

 Il y a trois données:
 
 \smallskip
  a) $k$  est un corps de nombres de clôture algébrique $\kbar$; on note $\Gamma_k$ le groupe de Galois $ \Gal(\kbar/k)$.
  
 b)  $L$ est un ensemble de nombres premiers.
  
  c) Pour tout $\l \in L$, $G_\l$ est un groupe de Lie $\l$-adique localement compact\footnote{On ne perdrait rien si l'on supposait que les $G_\l$ sont
compacts, vu que l'on peut supposer que les  $\r_\l$  sont surjectifs.}, et $\r_\l : \Gamma_k \to G_\l$ est un homomorphisme continu.

\smallskip

On fait deux sortes d'hypothèses:
  
  \smallskip
  
\n 2.1.   On suppose que la famille des  $\r_\l(\Gamma_k)$ est {\it bornée} , i.e. qu'elle satisfait
    à la condition suivante:
    
 \n   (B) {\it Il existe un entier  $n \geqslant 0$  tel que, pour tout  $\l \in L$, $\r_\l(\Gamma_k)$ soit isomorphe à un sous-quotient de $\GL_n(\Z_\l)$.}
    
    [Rappelons qu'un ``sous-quotient" d'un groupe  $A$ est un quotient d'un sous-groupe de $A$. Bien sûr, il s'agit ici de sous-groupes fermés.]
    
      Les cas particuliers les plus intéressants sont ceux où  $G_\l =
      \GL_{n_\l}(\Z_\l)$, ou $ G_\l = \GL_{n_\l}(\F_\l)$, avec des $n_\l$ bornés (par exemple constants).
      
      \smallskip
      
 \n     2.2. On fait une hypothèse du genre `` {\it semi-stabilité} '' sur la famille des $\r_\l$. 
      Pour l'énoncer, notons $V_k$ l'ensemble des places non archimédiennes de $k$. Si $v \in V_k$, notons  $k_v$ le complété de  $k$
      en  $v$,  notons $p_v$ la caractéristique résiduelle de  $v$ et choisissons un prolongement $\vbar$ de  $v$  à  $\kbar$. Notons $I_{\vbar}$ le groupe d'inertie correspondant à $\vbar$;
      c'est un sous-groupe fermé de $\Gamma_k$; à conjugaison près, il ne
      dépend que de $v$.

            Avec ces notations, l'hypothèse dont on a besoin
            s'énonce de la manière suivante:
      
   \n   (ST) {\it Il existe un sous-ensemble fini $S$ de $V_k$ tel que}:
         
  (ST1)    {\it  Si $v \notin S$ et $\l \neq p_v$, alors $\r_\l(I_{\vbar}) = 1$, i.e. $\r_\l$ est non ramifié en $v$.}
 
  (ST2) {\it Si $v \in S$ et $\l \neq p_v$, alors $\r_\l(I_{\vbar})$ est un pro-$\l$-groupe.}
  
\n  [Noter que l'on ne fait aucune hypothèse sur les $\r_\l(I_{\vbar})$ lorsque $\l=p_v$.]

\smallskip
  
   Lorsque $G_\l = \GL_n(\Z_\l)$, la condition (ST2) est moins restrictive
   que la condition habituelle de semi-stabilité, où l'on exige que $\r_\l(I_{\vbar})$
   soit formé d'éléments unipotents.
   
     Il est commode d'introduire une notion analogue à la {\it potentielle semi-stabilité} :
     
 \n    (PST) {\it Il existe une extension finie de  $k$ pour laquelle} (ST)  {\it est
     satisfaite.}
     [Plus explicitement: il existe une sous-extension finie $k_1$ de $\kbar$ telle
     que la famille des $\r_\l|\Gamma_{k_1}$ satisfasse à (ST).]
     
     Noter que, dès que (ST) est satisfaite pour une extension
     $k_1$ de $k$, elle l'est aussi pour toute extension finie de $k$ contenant $k_1$.
     
   \smallskip
   
\n   2.3. Le théorème que nous avons en vue dit que les propriétés
(B) et (PST) entraînent la propriété (PR) du \S1. Autrement dit:
   
   \smallskip
   
    \n {\bf  Théorème 1.}  {\it Si la famille des $\r_\l(\Gamma_k)$ est bornée au sens de $\rm(B),$ et si la condition  $\rm(PST)$ est satisfaite, il existe une extension
    finie de $k$ sur laquelle les $\r_\l$ sont indépendants.}
    
   \smallskip
    
    On peut reformuler cet énoncé en termes d'extensions de $k$: notons  $N_\l $ le noyau de $\r_\l$
    et $k_\l$ le sous-corps de $\kbar$ fixé par $N_l$; posons $N'_\l = \bigcap_{\l' \neq \l} N_{\l'}$ et notons $k'_\l$ le corps fixé par $N'_\l$, autrement dit le corps engendré par les $k_{\l'}$ avec $\l' \neq \l$. Le corps $k^{\ind} = \bigcap k'_\l$  correspond, par la théorie de Galois, au plus petit sous-groupe fermé
    de $\Gamma_k$ contenant les $N'_\l$. Avec ces notations, le théorème 1 est équivalent à:

    \smallskip
    
     \n {\bf  Théorème 1$'$.} {\it Si les conditions} (B) {\it et} (PST) {\it sont satisfaites,
     le corps $k^{\ind} = \bigcap_\l   k'_\l$ défini ci-dessus est une extension
     finie de $k$.}
     
     De plus, $k^{\ind} $ est la plus petite extension de $k$ sur laquelle les $\r_\l$
     sont indépendants; on peut l'appeler le ``corps d'indépendance" des $\r_\l$.
    
    \smallskip
    
    La démonstration des théorèmes 1 et 1$'$ sera donnée au \S8.
        
    \bigskip
    
  \n {\bf \S3. 
 Exemples et contre-exemples.}
 
   Dans chacun des exemples ci-dessous, l'ensemble  $L$  est l'ensemble de tous les
   nombres premiers, et $G_\l$ est isomorphe à $\GL_n(\Q_\l)$, avec  $n$
   fixe. Cette dernière hypothèse entraîne que le groupe  $\r_\l(\Gamma_k)$
   est isomorphe à un sous-groupe fermé de $\GL_n(\Z_\l)$, de sorte que
   la condition (B)  est satisfaite.
   \smallskip

\n 3.1. {\it Variétés abéliennes et quasi-abéliennes.}    Si $A$ est une variété abélienne de dimension  $d$ sur $k$, les modules de Tate $T_\l(A)$ fournissent des représentations $\l$-adiques de dimension $2d$
de $\Gamma_k$ qui satisfont à (PST) en vertu du théorème de Grothendieck et Mumford sur la semi-stabilité des modèles de Néron ([SGA 7 I, exposé IX], voir aussi [BLR 90, \S7.4]).

D'après le théorème 1, ces représentations sont presque indépendantes:
on retrouve ainsi un résultat démontré un peu différemment dans [Se 86].  Noter qu'ici les corps $k_\l$  ont une interprétation simple: $k_\l$ est le corps de rationalité des points de $A(\kbar)$ d'ordre une puissance de $\l$, et $k'_\l$ est le corps de rationalité des points de $A(\kbar)$ d'ordre fini premier à $\l$.

  Ces résultats s'appliquent aussi au cas des schémas en groupes
  quasi-abéliens; ce cas a été utilisé par Hrushovski, cf. [Bo 00].%\footnote{C'est d'ailleurs l'exposé Bourbaki de E. Bouscaren sur les travaux de Hrushovski qui m'a convaincu que je devais publier la démonstration du théorème 1.}

  \smallskip

\n 3.2. {\it Cohomologie $\l$-adique.}  Plus généralement, si  $X$  est un schéma séparé de type fini sur  $k$, la condition (PST) est satisfaite par les représentations $\l$-adiques associées aux {\it groupes de cohomologie} à support propre  $H^i_c(\Xbar,\Q_\l)$, ainsi que par les groupes de cohomologie $H^i(\Xbar,\Q_\l)$ à support quelconque. En effet:

  a) La condition  (ST1)  est satisfaite d'après les théorèmes d'existence de ``stratifications" dus à  N. Katz et G. Laumon [KL 86, th.3.1.2 et th.3.3.2].
  
  b) Si $S$ est choisi comme dans (ST1), il résulte d'un théorème de
Berthelot [Be 96, prop.6.3.2] que, pour tout  $v \in S$, il existe un sous-groupe ouvert normal $U_{\vbar}$       de $I_{\vbar}$ qui opère de façon unipotente
sur les  $H^i_c(\Xbar,\Q_\l)$  et les $H^i(\Xbar,\Q_\l)$, pourvu que  $\l \neq p_v$.
[La démonstration de Berthelot est basée sur  la théorie des altérations de de Jong, cf. [Jo 96].] Choisissons une extension galoisienne  $k'_v$ du corps local $k_v$ telle que    $\Gamma_{k'_v} \cap I_{\vbar} \subset U_{\vbar}$ . Un argument d'approximation bien connu montre qu'il existe une extension galoisienne finie $k_1$ de $k$ dont les complétés locaux aux places au-dessus de  $S$  contiennent les  $k'_v$. On a alors $\Gamma_{k_1} \cap I_{\vbar} \subset U_{\vbar}$ pour tout $ v \in S$, ce qui montre que la condition
(ST2)  est satisfaite sur  $k_1$.
  
 \smallskip

 \n {\it Problème} (cf. [Se 91, 10.1?]). Au lieu de supposer, comme nous venons de le faire, que $k$ est
 un corps de nombres, supposons seulement que $k$ est une extension de type fini de $\Q$. Comme ci-dessus, soit  $X$  un schéma séparé de type fini sur  $k$. Est-il encore vrai que les 
 représentations $\l$-adiques de $\Gamma_k$ fournies par les  $H^i_c(\Xbar,\Q_\l)$ et
 les $H^i(\Xbar,\Q_\l)$ sont presque indépendantes ? 
     
      \smallskip
    
 \n   3.3. {\it Mariage ``carpe-lapin"}. On peut partir de deux familles de $\r_\l$ satisfaisant aux hypothèses (B) et (PST),
    et pour chaque  $\l$  choisir au hasard l'un des deux $\r_\l$; on obtient encore une famille presque indépendante. Exemple: pour $\l \equiv 1 \pmod 4$ prendre la représentation $\l$-adique associée à la fonction de
    Ramanujan, et pour les autres $\l$ la représentation $\l$-adique associée
    à la courbe elliptique d'équation $y^2-y=x^3-x^2$.
    
    \smallskip
    \n 3.4. {\it Exemple montrant que la condition (PST) ne peut pas être entièrement supprimée.} 
    Soit $k = \Q$. Choisissons un nombre premier $p > 2$, ainsi qu'une suite infinie  $\l_1 < \l_2 < ... $ de nombres premiers tels
    que $\l_i   \equiv 1 \pmod {p^i}$. Soit $L = \{\l_1, \l_2, ...\}$. Soit $\r_{\l_i}: \Gamma_k  \to  \Z_{\l_i}^\times$ 
    un homomorphisme non ramifié en dehors de $p$ dont l'image est 
    cyclique d'ordre $p^i$. La famille des $\r_{\l_i}$  satisfait à la condition (B)
    avec $n=1$ et à la condition (ST1) avec $S = \{p\}$;  elle ne possède cependant pas la propriété (PR) car son corps d'indépendance est l'unique $\Z_p$-extension de $\Q$, qui est de degré infini sur $\Q$.

    \bigskip
    
     \n {\bf \S4. Un théorème de finitude sur les corps de nombres.}
     
       Soit $d$ un entier $> 0$, et soit $G$ un groupe fini. Considérons la condition:
       
       \smallskip
       
       \n (Jor$_d$) {\it Il existe un sous-groupe abélien normal  $A$ de  $G$  tel que} $(G:A) \leqslant d$.
       
       \smallskip
       
         \n {\bf  Théorème 2.} {\it Pour tout  $d > 0$ il n'existe qu'un nombre
         fini de sous-extensions galoisiennes  $K/k$  de $\kbar/k$  qui sont partout non ramifiées et dont le groupe de Galois a la propriété}   \rm({Jor}$_d)$  {\it  ci-dessus.}  
         
         \smallskip
         
 \n        {\it Démonstration.} On sait (Hermite-Minkowski) qu' il n'existe qu'un nombre fini de sous-extensions de $\kbar$ de degré $\leqslant d$ qui soient partout non ramifiées (cela provient de ce que leurs
         discriminants sont bornés en valeur absolue, cf. par exemple [Se 81, §1.4]). On peut donc trouver 
         une sous-extension finie $k_1$ de $\kbar$ contenant toutes ces extensions. Soit $k_2$ la plus grande extension abélienne non ramifiée
         de  $k_1$ contenue dans $\kbar$; d'après la théorie du corps de
         classes, $k_2$  est une extension finie de  $k_1$, donc aussi
         de $k$. Soit maintenant  $K/k$  une extension galoisienne dont le
         groupe de Galois $G$ a la propriété de l'énoncé, et soit  $K'$
         le sous-corps de  $K$  fixé par un sous-groupe abélien normal $A$ d'indice $\leqslant d$. On a $[K':k] \leqslant (G:A) \leqslant d$ et $K'$ est
         non ramifiée sur  $k$. Cela montre que  $K'$  est contenu dans  $k_1$. Comme  $K/K'$ est abélienne et non ramifiée, il en est de même
         de $K.k_1/k_1$ et cela entraîne que  $K.k_1$ est contenu dans
         $k_2$, d'où $K \subset k_2$, ce qui prouve
         la finitude cherchée.   
         
        [Ce théorème utilise deux des propriétés les plus importantes des
        corps de nombres:
        
         a) finitude des extensions de $\Q$  de degré et 
        discriminant bornés\footnote{En fait discriminant borné entraîne degré borné, mais cela ne joue aucun rôle ici.}; 
        
        b) finitude des extensions abéliennes non ramifiées (corps de classes).]

\bigskip

     \n {\bf \S5. Groupes linéaires d'ordre premier à la caractéristique.}
     
     \smallskip

\n 5.1. {\it Le théorème de Jordan classique.}

  Sous sa forme originelle ([Jo 78]), ce théorème s'énonce comme suit:
  
  \smallskip
  
    \n {\bf  Théorème 3.} {\it Pour tout entier $n \geqslant 0$ il existe un entier  $d = d(n)$ tel que tout sous-groupe fini de $\GL_n(\C)$ ait la propriété    \rm({Jor}$_d) \ du \ \S4.$}
    
    [Autrement dit, un sous-groupe fini de $\GL_n(\C)$ ne peut être ``gros"
    que s'il contient un gros sous-groupe abélien.]
    
    \smallskip
    
    On trouvera dans [Fr 11] une démonstration simple de ce résultat.
    Cette démonstration donne une valeur de  $d(n)$  telle que
       $$ d(n) \leqslant (\sqrt{8n}+1)^{2n^2} .$$
       On connaît maintenant la valeur optimale de  $d(n)$, qui est bien
       inférieure à celle-là; ainsi, pour $n \geqslant 71$, on a $d(n) = (n+1)!$, d'après M.J. Collins [Co 07], améliorant des résultats de B. Weisfeiler et de W. Feit \footnote{Les démonstrations de Weisfeiler, Feit et Collins dépendent de la
       classification des groupes finis simples.}. Nous n'en aurons pas besoin. Dans ce qui suit, nous noterons  $d(n)$ n'importe quel entier  $d$ pour lequel le théorème $3$
    est valable.
    
    \smallskip
    
    \n 5.2. {\it Le théorème de Jordan sur un corps quelconque.}
  \smallskip
  
  \n {\bf  Théorème 3$'$.} {\it Soient  $n$ un entier $
  \geqslant 0$,  $F$  un corps, $H$  un sous-groupe fini de $\GL_n(F)$  et  $G$  un quotient de $H$. On suppose que $|G|$  est premier à la caractéristique de $F$ si celle-ci est $\neq 0$. Alors  $G$  a la propriété   \rm({Jor}$_{d(n)}) \ du \ \S4$}.
  
  \smallskip

  \n {\it Démonstration.} Elle se fait en trois étapes:
  \smallskip
  
  \n 5.2.1. Le cas où car($F)=0.$  On peut supposer  $F$  de type fini sur $\Q$, donc plongeable dans $\C$. Le théorème 3 montre alors que $H$ a la propriété  \rm({Jor}$_{d(n)})$ et il en est donc de même de $G$.
  
  \smallskip
  
  \n 5.2.2. Le cas où car($F) = p > 0$ , avec  $|H|$  premier à $p$. On peut
  supposer que  $F$  est parfait. Soit $W$ l'anneau des vecteurs de Witt à coefficients dans  $F$. On a
  un homomorphisme surjectif   $\GL_n(W) \to \GL_n(F)$. Comme  $|H|  $  
est premier à $p$, $H$ se relève en un sous-groupe de $\GL_n(W)$, 
et l'on applique 5.2.1 au corps des fractions de $W$.

\smallskip

\n 5.2.3. Le cas où car($F) = p > 0$. Soit $I$ le noyau de $H \to G$, et soit $P$ un $p$-Sylow de $I$; c'est aussi un $p$-Sylow de $H$, puisque $(H:I)$
est premier à $p$. Soit $N_H(P)$ le normalisateur
de $I$ dans $H$. On sait (Frattini) que $N_H(P) \to G$ est surjectif\footnote{L'argument dit ``de Frattini"  est le suivant: si $h\in H$, $hPh^{-1}$ est un $p$-Sylow de $I$, donc s'écrit  $xPx^{-1}$ avec  $x \in I$, d'où
$x^{-1}h \in N_H(P)$, ce qui montre que  $h$  appartient à $I.N_H(P)$.
On a donc bien $H = I.N_H(P).$}. D'autre part, la suite exacte
$$  1 \to P  \to N_H(P) \to N_H(P)/P \to 1$$
est scindée car les ordres de $P$ et de $N_H(P)/P$ sont premiers entre eux. Il existe donc un sous-groupe $H'$ de $N_H(P)$, d'ordre premier à $p$, tel que $N_H(P) = P.H'$. Or l'image de $P$ dans $G$  est triviale,
puisque $P$ est contenu dans $I$. On en déduit que $G$ est un quotient
de $H'$, et l'on conclut en appliquant 5.2.2 à $H'$.

\bigskip
 
   \n {\bf \S6. Groupes linéaires engendrés par des éléments d'ordre égal à la caractéristique.}
   
     Dans ce qui suit, $\l$ désigne un nombre premier $\geqslant 5$. 
     
   \smallskip
   
       \n 6.1. {\it  Les groupes simples finis de caractéristique $\l$}: {\it la famille $\Sigma_\l$. }
       
  \smallskip
    
     Rappelons comment on définit les groupes simples  ``du type de Lie"
     en caractéristique $\l \geqslant 5$ (pour les propriétés utilisées ici,
     voir par exemple [GLS 98, \S2.2] - noter que l'hypothèse $ \l \geqslant 5$
     élimine les cas particuliers exceptionnels que l'on rencontre en caractéristique 2 et 3, ainsi que les formes tordues à la Suzuki-Ree).
     
      On se donne un groupe algébrique lisse connexe $\underline{H}$ sur un corps fini $F$ dont l'ordre est une puissance de $\l$. On suppose
     que  $\underline{H}$ est géométriquement simple et simplement connexe, et l'on désigne par $\underline{H}^{adj}$ le quotient de $\underline{H}$ par son centre. L'image $H_F$ de
     l'homomorphisme  $\underline{H}(F) \to  \underline{H}^{adj}(F)$ est alors un groupe fini simple
     non abélien.
     
     \smallskip
     \n {\it Remarque.} On aurait aussi pu définir 
 $H_F$ comme le quotient de $\underline{H}(F)$ par son centre, ou bien comme le sous-groupe de $\underline{H}^{adj}(F)$ engendré par 
     les $\l$-Sylow de $\underline{H}^{adj}(F)$. L'équivalence de ces diverses définitions provient de ce que $\underline{H}(F)$ est engendré par ses
     éléments unipotents d'après un théorème de Steinberg [St 68, th.12.4].
     
     \smallskip
     
       Nous noterons $\Sigma_\l$ l'ensemble des classes d'isomorphisme
       de groupes finis simples qui sont, soit du type $H_F$ ci-dessus (pour  un $\underline H$ et un
       $F$ convenables\footnote{Il y a unicité: un groupe simple n'est isomorphe à $H_F$ que pour au plus un couple $(\underline H,F)$, à isomorphisme près.}), soit isomorphe au groupe cyclique $\Z/\l\Z$.

     \smallskip
     
      \n 6.2. {\it  Un lemme.}
      \smallskip
      
 \n {\bf  Lemme 1.} {\it Soit  $\underline{G}$ un groupe algébrique linéaire connexe sur $\F_\l$ et soit  $G = \underline{G}(\F_\l)$ le groupe de ses points rationnels. Tout quotient simple d'une suite de Jordan-Hölder de $G$
 appartient\footnote{Dans ce qui suit, on dit qu'un groupe simple ``appartient" à  $\Sigma_\l$ lorsqu'il est isomorphe à un élément
 de $\Sigma_\l$.} à $\Sigma_\l$ ou est cyclique d'ordre  $\neq \l$.}
     
 \n    {\it Démonstration}. Un argument de dévissage permet de supposer
     que  $\underline{G}$ est, soit un groupe unipotent, soit un tore, soit un groupe semi-simple. Les deux premiers cas sont immédiats. On peut donc supposer que  $\underline{G}$ est semi-simple.  Soit $\widetilde{\underline{G}}$ le revêtement universel de $\underline{G}$ et soit $\underline{G}^{adj}$ son groupe adjoint. Soient $\widetilde{G}$ et $G^{adj}$ les groupes de points $\F_\l$-rationnels de ces groupes algébriques. On a des homomorphismes
     naturels $$\widetilde{G} \ \to \ G \ \to  \ G^{adj}.$$
     Comme $\widetilde{\underline{G}}$ est simplement connexe, c'est un produit de groupes du type $R_{F/{\scriptsize\F}_\l} \underline{H}$, où $\underline{H}$ et $F$ sont comme dans 6.1 ci-dessus, et
     le symbole $R_{F/{\scriptsize\F}_\l}$ désigne le foncteur ``restriction des scalaires" à la Weil
     (celui que Grothendieck note $\prod_{F/{\scriptsize\F}_\l}$), cf. par exemple [KMRT 98, th.26.8]. On a donc $\widetilde{G}= \prod \underline{H}(F)$.   Les homomorphismes
     $$ \widetilde{G} \ \to \ G \ \to  \ G^{adj}$$
     
   \n  ont des noyaux et conoyaux qui sont commutatifs d'ordre premier à $\l$.
     De plus, l'image de $ \widetilde{G} $ dans $G^{adj}$  est un produit de
     groupes simples appartenant à  $\Sigma_\l$. Le lemme en résulte.
     
     \smallskip

\medskip

  \n 6.3. {\it Un théorème de Nori.}
  \
  \smallskip
  
    \n {\bf  Théorème 4.} {\it  Pour tout $n \geqslant 0$, il existe un entier
    $c(n)$ tel que, si $\l > c(n)$, tout sous-quotient fini simple de $\GL_n(\Z_\l)$ d'ordre divisible par $\l$ appartient à $\Sigma_\l$.}
     
     \smallskip
     
     \n{\it Démonstration.} Prenons $c(n)= \sup{(3,c_2(n))}$, où $c_2(n)$ a les propriétés énoncées dans  [No 87, Theorem B]. Nous allons voir que cet entier convient. 
     
     Supposons que  $\l>c(n)$ et soit $H$ un sous-quotient fini simple de $\GL_n(\Z_\l)$ d'ordre divisible par $\l$. Comme $H$ est simple, cette dernière propriété entraîne que $H$ est engendré par ses
     $\l$-Sylow. 
     
    L'homomorphisme naturel 
    $\GL_n(\Z_\l)  \ \to \  \GL_n(\F_\l)$ 
    est surjectif, et son noyau est un pro-$\l$-groupe. Il en résulte que $H$ est, soit cyclique d'ordre $\l$, soit isomorphe à un sous-quotient de $\GL_n(\F_\l)$. Dans le 
    premier cas, $H$ appartient à  $\Sigma_\l$. Dans le second cas, on a $H = G/I$, avec $G \subset \GL_n(\F_\l)$ et $I$ normal dans $G$; on peut évidemment supposer que $G$ est engendré par ses $\l$-Sylow. D'après  [No 87, Th.B],
     il existe un $\F_\l$-sous-groupe algébrique connexe $\underline{G}$ de $\GL_n$ tel que $G$ soit contenu dans $\underline{G}(\F_\l)$ et soit engendré par les $\l$-Sylow de ce groupe\footnote{  La définition de $\underline{G}$ donnée par Nori est très simple: c'est le plus petit sous-groupe
     algébrique de $\GL_n$ contenant les groupes à 1 paramètre $t \mapsto u^t$, où
    $u$ parcourt les éléments d'ordre $\l$ de $G$. Dans la terminologie
    de [Se 94, \S4], c'est le {\it saturé} de $G$.}. Le groupe $H$ est un quotient d'une suite de Jordan-Hölder de $G$, donc aussi de $\underline{G}(\F_\l)$.
     D'après le lemme 1, ceci entraîne que $H$ est, soit cyclique d'ordre
     premier à $\l$ (ce qui est exclu), soit isomorphe à un élément de  $\Sigma_\l$. D'où le théorème.
     
  \medskip
 
   \n 6.4. {\it Un théorème d'Artin.}
  
     Le  résultat suivant est essentiellement dû à E. Artin ([Ar 55], complété par [KLST 90]):
       
       \smallskip
         \n {\bf  Théorème 5.} {\it Si  $\l' $ est premier  $\geqslant 5$ et distinct de $\l$, 
         on a} \ \  $\Sigma_\l \cap \Sigma_{\l'} = \varnothing$.
         
          \smallskip
          
          La démonstration donne même un résultat plus fort: si $G$ appartient à $\Sigma_\l$ et $G'$ appartient à $ \Sigma_{\l'}$, leurs ordres $|G|$ et  $|G'|$ sont distincts.          
          \smallskip
          
          \n {\it Exemples}. Pour $\l = 5$, les ordres des éléments de $\Sigma_\l$,
          rangés par taille croissante, sont 
          $\{5, 60, 7800, 126000, 372000, 976500, ...\}$.
          
          \n Pour $\l=7$, ce sont $\{7, 168, 58800, 1876896, 5663616, 20176632, ...\}.$
          
           \bigskip
           
             \n {\bf \S7. Deux critères d'indépendance.}
             
             \smallskip 
 
 \n 7.1. {\it Un critère élémentaire.}
 
 \smallskip
 
   Revenons aux notations du \S1, et soit $\r_i : \G \to G_i, \ i\in I$, une famille d'homomorphismes, les groupes $\G$ et $G_i$ étant des groupes
   profinis, et les  $\r_i$ étant continus.
   
   \smallskip
   
    \n {\bf  Lemme 2.} {\it Supposons que les groupes $\r_i(\G) \subset G_i$
    aient la propriété suivante}:
    
    (D) {\it Si $i \neq j$, aucun quotient fini simple de $\r_i(\G)$ n'est isomorphe à un quotient de $\r_j(\G)$.
     
\n     Alors les $\r_i$ sont indépendants.}
     
     \smallskip
     
     \n {\it Démonstration.} On peut évidemment supposer que  les $\r_i$ sont surjectifs, i.e.  $G_i = \r_i(\G)$ pour tout  $i$.
     
     Considérons d'abord le cas où $I$ est un ensemble à
     deux éléments, par exemple $I = \{1,2\}$. Si $\r: \G \to G_1 \times G_2$
     n'est pas surjectif, le classique lemme de Goursat montre qu'il existe un groupe
     profini non trivial $A$ et des homomorphismes surjectifs $f_i:G_i \to A$
     tels que $f_1 \circ \r_1 = f_2 \circ \r_2$. Comme  $A$  est non trivial, il a
     un quotient qui est un groupe simple fini, et ce groupe est quotient à la
     fois de $G_1$ et de $G_2$, contrairement à l'hypothèse (D).
     
     Le cas où $I$ est fini se déduit par récurrence sur $|I|$ du cas où  $|I| = 2$, et le cas où $|I|$ est infini se déduit par passage à la limite du cas où $|I|$ est fini.
 
 \medskip
 
  \n 7.2. {\it Un autre critère}.
 
   Soit $\G$ un groupe profini et soit $L$ un ensemble de nombres premiers. Pour tout $\l \in L$, soit $\r_\l : \G \to G_\l$ un homomorphisme continu de $\G$ dans un groupe de Lie $\l$-adique compact $G_\l$.
   
   \smallskip
   
    \n {\bf  Lemme 3.} {\it Supposons qu'il existe une partie finie $I$ de $L$
    telle que la famille $(\r_\l)_{\l \in L - I}$ ait la propriété $\rm{(PR)}$ du $\S1.$
    Alors il en est de même de la famille $(\r_\l)_{\l \in L}$.}
    
     (Autrement dit, pour prouver (PR), on a le droit de supprimer un nombre fini
     d'éléments de $L$.)
     
     \smallskip
     
         \n {\it Démonstration.} On peut supposer que $I$ est réduit à un seul élément, que l'on notera $p$: le cas général en résultera par récurrence sur $|I|$. 
         Quitte à remplacer $\G$ par un sous-groupe ouvert, on peut  supposer que les $\r_\l$ sont indépendants pour $\l \neq p$; on peut aussi supposer
         que tous  les $\r_\l$ sont surjectifs. Nous allons alors démontrer un peu mieux
         que (PR), à savoir:
         
         \smallskip
         
         (*) {\it La famille des $\r_\l$ possède la propriété $\rm(RO)$ du} \S1.
         
         \smallskip
         
   \n      Autrement dit, l'image de  $\G$  par l'homomorphisme
         
         $$  \r = (\r_\l): \G  \  \to  \  G_p \times \prod_{\l \neq p}  G_\l $$
         
       \n  est ouverte dans $\prod_\l G_\l$.
         
         Les deux projections $\r(\G) \to G_p$ et $\r(\G) \to \prod_{\l \neq p} G_\l $ sont surjectives par hypothèse. On se trouve donc dans la situation du lemme de Goursat. Autrement dit,
          si l'on identifie  $G_p$  au facteur  $G_p \times 1$ de $G_p \times \prod_{\l \neq p} G_\l $, le groupe quotient  $C = G_p/(\r(\G) \cap G_p)$ est un quotient de $\prod_{\l \neq p} G_\l $. Dire que
         $\r(\G)$ est ouvert équivaut à dire que $C$ est {\it fini}.
         C'est ce que nous allons démontrer.
         
           Observons d'abord que $C$ est un groupe de Lie $p$-adique compact (puisque c'est un quotient de $G_p$); il contient donc un sous-groupe ouvert normal  $U$ qui est un pro-$p$-groupe sans torsion (cf. par exemple [Se 65, II, \S IV.9, th.5], [Bo 72, Chap.III, \S7] ou [DSMS 99, th.8.32]). \ Si  \  $J$  \ est        
           une partie finie  de $L \ \se \ \{p\}$,  notons $C_J$ l'image de l'homomorphisme   $$\prod_{\l \in J} G_\l  \  \ \to \prod_{\l \neq p} G_\l  \  \to \  C. $$
         Les $p$-Sylow des $G_\l$ sont finis si  $\l \in J$; il en est donc
         de même de ceux de $C_J$. Comme  $U$  est sans torsion, cela
         montre que $U \cap C_J = 1$; d'où $|C_J| \leqslant (C:U)$. Cela donne 
         une borne uniforme pour l'ordre de $C_J$, ce qui entraîne qu'il existe
         un $C_J$ qui contient tous les autres. Mais la réunion des $C_J$ est dense dans $C$.
          D'où le fait que  $C$  est fini.
     
\bigskip

 \n {\bf \S8. Démonstration du théorème 1.}
 
 Revenons à la situation du théorème 1, relative à un homomorphisme
 $$\r = (\r_\l): \Gamma_k \to \prod_{\l \in L} G_\l$$ satisfaisant aux conditions
 (B) et (NST). Pour prouver que $\r$ a la propriété (PR), nous procéderons en plusieurs étapes.
 
 \smallskip
 
 \n 8.1. {\it Réductions.} Quitte à remplacer $k$ par une extension finie, on peut supposer que la condition de semi-stabilité (ST) est satisfaite. On peut aussi 
 supposer que les $\r_\l$ sont surjectifs. D'après (B),
 on peut choisir un entier $n \geqslant 0$ tel que, pour tout  $\l \in L$, le groupe $G_\l$ soit un sous-quotient de $\GL_n(\Z_\l)$. D'après le
 lemme 3, on peut aussi supposer que tous les $\l \in L$ sont $> \sup(3,c(n)) $  où $c(n)$ a la propriété énoncée dans le théorème 4. Pour la même raison, on peut
 aussi supposer que l'on a $\l \neq p_v$  pour toute place  $v$  de l'ensemble fini $S$ intervenant dans (ST). 
 
 \smallskip
 
\n 8.2. {\it Les groupes $A_\l$.} Si $\l \in L$, notons $\Gamma_{k,\l}$ le plus petit sous-groupe normal fermé de $\Gamma_k$ contenant les groupes d'inertie $I_{\vbar}$ correspondant aux places $v$ telles que $p_v = \l$.  D'après (ST1), on a $\r_{\l'}(\Gamma_{k,\l}) = 1$ pour tout $\l' \neq \l$. 
L'image de $\Gamma_{k,\l}$ par $\r: \Gamma_k \to \prod G_\l$ est donc contenue dans le $\l$-ième facteur de $\prod G_\l$. Notons $A_\l$ cette
image; c'est un sous-groupe fermé normal de $G_\l$. Le plus petit sous-groupe fermé de $\prod G_\l$ contenant tous les
$A_\l$ n'est autre que le produit $\prod A_\l$. En particulier, on a:

\smallskip

 \n {\bf  Lemme 4.}  {\it Le sous-groupe $\r(\Gamma_k)$ de $\prod G_\l$ contient $\prod A_\l$.}
 
 \smallskip
 
 \n 8.3. {\it Les groupes $G^+_\l$.} Si $\l \in L$, notons
$G^+_\l$ le sous-groupe de $G_\l$ engendré par ses $\l$-Sylow;
  c'est un sous-groupe ouvert normal de $G_\l$.  Posons $H_\l =  G_l/G_\l^+ . A_\l$;
  c'est un groupe fini d'ordre premier à $\l$.
 
  \smallskip
  
   \n {\bf  Lemme 5.}   a) {\it L'homomorphisme $\Gamma_k \to G_\l \to   H_\l$ est partout non ramifié.}
  
\n  b) {\it Le groupe $H_\l$ jouit de la propriété} Jor$_{d(n)}$ {\it des} \S\S 4-5.
  
  \smallskip
  
  \n {\it Démonstration.} Soit  $v \in V_k$, et soit $\vbar$ une place de $\kbar$
  prolongeant $v$. Si $p_v = \l$, on a $\r_\l(I_{\vbar}) \subset A_\l$  par définition de  $A_\l$; l'image de $I_{\vbar}$ dans $H_\l $ est donc triviale. Si $p_v \neq \l$,
  le groupe  $\r_\l(I_{\vbar})$  est un pro-$\l$-groupe d'après (ST); il est donc contenu dans $G_\l^+$  et son image dans $H_\l $ est triviale. Cela démontre a).
  
    Quant à b), il résulte du fait que l'ordre de $H_\l$ est premier à $\l$, ce qui
    permet de lui appliquer le théorème 3$'$.
    
    \smallskip
    
 \n   8.4. {\it Changement de corps.} D'après le lemme 5, les homomorphismes
    $\Gamma_k \to H_\l$ sont non ramifiés. Comme les $H_\l$ ont la propriété Jor$_{d(n)}$, on peut appliquer le théorème 2. On en déduit qu'il existe une extension finie non ramifiée  $k'$  de  $k$  telle que, pour tout $\l\in L$, l'image de $\r_\l(\Gamma_{k'})$ \
    dans  $H_\l$  soit triviale. Choisissons une telle extension.
    On a alors $\r_l(\Gamma_{k'}) \subset G_\l^+ . A_\l$ pour tout $\l$.
    Nous allons maintenant prendre  $k'$  comme corps de base; nous
    poserons  $G'_\l = \r_\l(\Gamma_{k'})$,
     et nous noterons  $G'^+_\l$ et
    $A'_\l$ les groupes correspondant à $G_\l^+$ et à $A_\l$; par exemple, $G'^+_\l$ est le sous-groupe de $G'_\l$ engendré par les $\l$-Sylow de $G'_\l$.
    
    \smallskip
    
    \n {\bf  Lemme 6.}  {\it Si $\l > [k':k]$, on a \ $G'^+_\l
 = G_\l^+$, \ $A'_\l = A_\l$ \ et \ $G'_\l = G'^+_\l.A'_\l.$ }
     
     \smallskip
     
    \n {\it Démonstration.}  L'hypothèse faite sur $\l$ entraîne que l'indice de $G'_\l$ dans $G_\l$ est $< \l$, d'où le fait que tout $\l$-Sylow de  $G_\l$  est contenu dans $G'_\l$, ce qui entraîne $G'^+_\l = G_\l^+$. L'égalité  $A'_\l = A_\l$ résulte de ce que les groupes d'inertie  $I_{\vbar}$  sont les mêmes
    pour  $k'$ et pour $k$, puisque  $k'$ est non ramifié sur  $k$. Enfin,
    l'égalité $G'_\l = G'^+_\l.A'_\l$ résulte de ce que  $G'_\l = \r_\l(\Gamma_{k'})$ est contenu dans $G_\l^+ . A_\l$. 
    
    \smallskip
    
\n 8.5. {\it Fin de la démonstration.}  D'après le lemme 3, on peut supposer que l'on a $\l > [k':k]$ pour tout $\l \in L$. Le lemme 6 montre que l'on a alors
$G'_\l = G'^+_\l.A'_\l$ pour tout  $\l$. D'après le théorème 4, tout quotient simple de  $G'^+_\l$  appartient à l'ensemble $\Sigma_\l$ défini au \no 6.1.
Il en est donc de même des quotients simples de $G'_\l/A'_\l$. Comme les
$\Sigma_\l$ sont deux à deux disjoints (théorème 5), on peut appliquer le lemme 2 à la famille des homomorphismes  $\Gamma_{k'} \to G'_\l/A'_\l$.
On en conclut que l'homomorphisme $\Gamma_{k'} \to \prod G'_\l/A'_\l$ est surjectif. Si l'on pose $X' = \r(\Gamma_{k'})$ et $A' = \prod A'_\l$, cela revient à dire que $X'.A' = \prod G'_\l$. Mais le lemme 4, appliqué au corps $k'$, montre que $X'$ contient $A'$. On a donc $X' = \prod G'_\l$, ce qui achève la démonstration.

\begin{center}
       {\bf Références}
\end{center}

[Ar 55] E. Artin, {\it The orders of the classical simple groups}, Comm. Pure
and Applied Math. $\bf 8$ (1955), 455-472 (= C.P., \no 33).

[Be 96] P. Berthelot, {\it Altération des variétés algébriques} ({\it d'après A.J. de Jong}), Sém. Bourbaki 1995/1996, exposé $\bf{815}$; Astérisque $\bf{241}$, SMF, 1997, 273-311.

[BLM 90] S. Bosch, W. Lütkebohmert \& M. Raynaud, {\it Néron Models}, Ergebnisse der
Math. (3) $\bf 21$, Springer-Verlag, 1990.

[Bo 72] N. Bourbaki, {\it Groupes et Algèbres de Lie, Chap.II et Chap.III}, Hermann, Paris, 1972.

[Bo 00] E. Bouscaren, {\it Théorie des modèles et conjecture de Manin-Mumford} ({\it d'après Ehud Hrushovski}), Sém. Bourbaki 1999/2000, exposé $\bf{870}$; Astérisque $\bf{276}$, SMF, 2002, 137-159.

%[DJ 86] A.J. de Jong, {\it Smoothness, semi-stability and alterations}, Publ. Math. IHES  $\bf{83}$ (1986), 51-93.

[Co 07] M.J. Collins, {\it On Jordan's theorem for complex linear groups}, J.
of Group Theory {\bf 10} (2007), 411-423.

[DSMS 99] J.D. Dixon, M.P.F. du Sautoy, A. Mann \& D. Segal, {\it Analytic pro-$p$-groups}, Second edition, revised and enlarged by Marcus du Sautoy \& Dan Segal, Cambridge Univ. Press, Cambridge, 1999.

[Fr 11] F.G. Frobenius, {\it \"{U}ber den von L. Bieberbach gefundenen Beweis eines Satzes von C. Jordan}, Sitz. Königlich Preuss. Akad. Wiss. Berlin (1911), 241-248 (= Ges. Abh.III, 493-500).

[GLS 98] D. Gorenstein, R. Lyons \& R. Solomon, {\it The Classification of the Finite Simple Groups, Number 3}, Math. Surveys and Monographs $\bf {40}$-3, AMS, 1998.

%[Gr 09] A. Greicius, {\it Elliptic curves with surjective adelic Galois representations}, arXiv:math/0901.2513v1 [math.NT.].

%[Gu 99] R.M. Guralnick, {\it Small representations are completely reducible}, J. Algebra $\bf220$ (1999), 531-541.

[Il 10] L. Illusie, {\it Constructibilité générique et uniformité en $\l$}, Orsay, 2010, non publié.
 
 [Jo 96] A.J. de Jong, {\it Smoothness, semi-stability and alterations}, Publ. Math. IHES 83 (1996), 51-93.
 
 [Jo 78] C. Jordan, {\it Mémoire sur les équations différentielles linéaires
 à intégrale algébrique}, J. Crelle $\bf84$ (1878), 89-215 (= Oe.II, 13-140).
 
 [KL 86] N.M. Katz \& G. Laumon, {\it Transformation de Fourier et majoration de sommes exponentielles}, Publ. Math. IHES 62 (1986), 361-418; {\it Erratum}, Publ. Math. IHES 69 (1989), 233.
 
 [KLST 90] W. Kimmerle, R. Lyons, R. Sandling \& D.N. Teague, {\it Composition factors from the group ring and Artin's theorem on orders of simple groups}, Proc. LMS $\bf60$ (1990), 89-122.
 
 [KMRT 98] M-A. Knus, A. Merkurjev, M. Rost \& J-P. Tignol, {\it The Book of Involutions}, AMS Colloquium Publ. 44, 1998.

[No 87] M.V. Nori, {\it On subgroups of $\GL_n(\F_p)$}, Invent. math. $\bf88$ (1987), 257-275.

%[Se 64] J-P. Serre, {\it Cohomologie Galoisienne}, Lect. Notes in Math. $\bf{5}$, Springer-Verlag, 1964; cinquième édition, révisée et complétée, 1994.

[Se 65] J-P. Serre, {\it Lie Algebras and Lie Groups}, Benjamin Publ., New York, 1965;  Lect. Notes in Math. $\bf1500$, Springer-Verlag, 1992; corrected fifth printing, 2006.

%[Se 72] J-P. Serre, {\it Propriétés galoisiennes des points d'ordre fini des courbes elliptiques}, Invent. math. $\bf15$ (1972), 259-331 (= Oe.III, \no 94).

[Se 81] ------ , {\it Quelques applications du théorème de densité de Chebotarev}, Publ. Math. I.H.E.S. $\bf54$ (1981), 123-201 (= Oe.III, \no 125).

[Se 86]  ------ , {\it Lettre à Ken Ribet du 7/3/1986} (= Oe.IV, \no 138).

[Se 91] ------ , {\it Propriétés conjecturales des groupes de Galois motiviques et des représentations $\l$-adiques}, Proc. Symp. Pure Math. $\bf55$, AMS, 1994, vol.I,
377-400 (= Oe.IV,  \no 161).

[Se 94]  ------ , {\it Sur la semi-simplicité des produits tensoriels de représentations de groupes}, Invent. math. $\bf116$ (1994), 513-530 (= Oe.IV, \no 164).

   [SGA 4] M. Artin, A. Grothendieck \& J-L. Verdier, {\it Théorie des Topos et Cohomologie  \'Etale des Schémas}, 3 vol., Lect. Notes in Math. $\bf269, 270, 305$, Springer-Verlag, 1972-1973.
   
   [SGA 7 I] A. Grothendieck, {\it Groupes de Monodromie en Géométrie Algébrique}, Lect. Notes in Math. $\bf{ 288}$, Springer-Verlag, 1972.

[St 68] R. Steinberg, {\it Endomorphisms of linear algebraic groups}, Memoirs AMS
$\bf80$ (1968) (= C.P., \no 23).

\bigskip

\n Collège de France, 3 rue d'Ulm, F-75005 Paris

\n serre@noos.fr

 \end{document}